\documentclass[12 pt]{article}
\usepackage{fancyhdr}

\newtheorem{Theorem}{Theorem}
\newtheorem{Lemma}{Lemma}

\begin{document}
\title{Functional central limit theorems for certain statistics in an infinite urn scheme}
\author{Mikhail Chebunin\thanks{E-mail: chebuninmikhail@gmail.com, Novosibirsk State University, Novosibirsk, Russia},
Artyom Kovalevskii \thanks{E-mail: kovalevskiii@gmail.com, Novosibirsk State Technical University, 
Novosibirsk State University, Novosibirsk, Russia}}
\date{}
\maketitle

\begin{abstract}
We investigate a specific infinite urn scheme first considered by Karlin (1967). 
We prove functional central limit theorems  for the total number of urns 
with at least $k$ balls for different $k$.
\end{abstract}


\section{Introduction}

Karlin (1967) studied an infinite urn scheme, that is, 
each of $n$ balls goes to urn $i\ge 1$ with probability $p_i>0$, $p_1+p_2+\ldots=1$, independently of other balls.
We assume $p_1\geq p_2 \geq \ldots $. 
Let $X_j$ be the box that the ball $j$ is thrown into,
and
\[
R^{*}_{n,k}=\sum\limits_{i=1}^{\infty} {\bf I} (\exists j_1<\ldots < j_k\leq n: \ X_{j_1}=\ldots=X_{j_k}=i)
\]
be the total number of urns with at least $k$ balls.
The number of nonempty urns is $R_n=R^{*}_{n,1}$. 
The total number of urns with exact $k$ balls is
$R_{n,k}=R^{*}_{n,k}-R^{*}_{n,k+1}$.
Let $ J_i (n) $ be the number of $n$ balls in urn $ i $.

Let (see Karlin (1967)) $ \Pi =\{\Pi(t),\ t \geq 0\} $ be a Poisson process with parameter 1.
This process does not depend on $\{X_j\}_{j\ge 1}$.
The Poissonized version of Karlin model assume the total number of $\Pi(n)$ balls.

According to well-known thinning property of Poisson flows,
stochastic processes
$\{J_i(\Pi(t))\stackrel{def}{=}\Pi_i(t), \ t\geq 0\}$  are Poisson with
 intensities $ p_i $ and are mutually independent for different $ i $'s. The definition implies that
\[
R^{*}_{\Pi(n),k}=\sum_{i=1}^{\infty} {\bf I}(\Pi_i(n)\geq k), \ \  R_{\Pi(n),k}
=\sum_{i=1}^{\infty} {\bf I}(\Pi_i(n)= k).
\]
Let $\alpha(x)=\max\{j|\ p_j\geq 1/x\}$
and  we assume
$\alpha(x)=x^{\theta} L(x)$, $0\le \theta\leq 1$, as in Karlin (1967). 
Here
$L(x)$ is a 
slowly varying function as $ x \to \infty $.
Let for $t\in [0,1], \ k\geq 1$ 
\[
Y_{n,k}^{*}(t) =\frac{R^{*}_{[nt],k}-{\bf E} R^{*}_{[nt],k}}{(\alpha(n))^{1/2}}, 
\ \ \ \ \ \ \ \ \ \  Z_{n,k}^{*}(t) =\frac{R^{*}_{\Pi (nt),k}-{\bf E} R^{*}_{\Pi (nt),k}}{(\alpha(n))^{1/2}},
\]
\[
 Y_{n,k}(t) =\frac{R_{[nt],k}-{\bf E} R_{[nt],k}}{(\alpha(n))^{1/2}},   \ \ \ \ \ \ \ \ \
K_{k,\theta}=\left\{
\begin{array}{ll}
-\Gamma(1-\theta), & k=0; \\
\theta \Gamma(k-\theta), & k>0,
\end{array}
\right.
\]

The goal of our paper is to extend the following two theorems from Karlin (1967).

\begin{Theorem} \ \ (Theorem 4 in Karlin (1967)).
Let \ $\theta\in(0,1]$. \ Then \ \ \ \ \ \ \ \ \ \ \ \ \ \ \ $ (R_n-{\bf E} R_n)/B_n^{1/2}$
converges weakly to standard normal distribution, where
\[
B_n=\left\{
\begin{array}{cr}
\Gamma(1-\theta)(2^{\theta}-1)n^{\theta}L(n), & \theta\in (0,1);\\
n \int\limits_0^{\infty}\frac{e^{-1/y}}{y}L(n y)dy \stackrel{def}{=}n L^{*}(n), & \theta=1.
\end{array}
\right.
\]
\end{Theorem}

Karlin (1967, Lemma 4) proved that the function $L^*(x)$ is
slowly varying as $x \to \infty$.

\begin{Theorem} \  (Theorem 5 in Karlin (1967)). 
Let $\theta\in(0,1)$, $r_1<\ldots <r_{\nu}$ 
be $\nu$ positive integers. Then random vector
$
\left(Y_{n,r_1}(1), 
\ldots,
Y_{n,r_\nu}(1) \right)
$
converges weakly to the multivariate normal distribution with zero expectation
 and covariances
\[
c_{r_i,r_j}=\left\{
\begin{array}{cr}
-\frac{\theta \Gamma(r_i+r_j-\theta)}{r_i! r_j!} 2^{\theta-r_i-r_j}, & i \neq j;\\
\frac{\theta}{\Gamma(r_i+1)}\left(
\Gamma(r_i-\theta)- 2^{-2r_i+\theta}\frac{\Gamma(2r_i-\theta)}{\Gamma(r_i+1)}\right), & i=j.
\end{array}
\right.
 \]
\end{Theorem}

Here we briefly mention some related results on this model.
Dutko (1989) generalized Theorem~1 by proving asymptotic normality of
$ R_n $ if ${\bf Var} R_n \to \infty$
as $n \to \infty$.  
This condition always holds if $\theta\in(0,1]$ but can hold too for $\theta=0$.
Gnedin, Hansen and Pitman (2007) focused on study of conditions for convergence 
${\bf Var} R_n \to \infty$. 
Barbour and Gnedin (2009) extended Theorem~2 on the case of $ \theta = 0 $ if
variances go to infinity. They  found conditions for convergence of
covariances to a limit
and identified four types of limiting behavior of variances.
Barbour (2009) proved theorems on approximation of the number of cells with $ k $ balls
by translated Poisson distribution.
Key (1992, 1996)  studied the limit behavior of statistics $R_{n,1}$.
Hwang and Janson (2008) 
proved local limit theorems for finite and infinite number of cells.
Zakrevskaya and Kovalevskii (2001)
proved consistency for one parametric family of an  estimator of $ \theta \in (0,1)$ which is an implicit function of $ R_n $.
Chebunin (2014) constructed an $ R_n $-based explicit parameter estimator for $\theta \in (0,1)$ and proved its consistency.
Durieu and Wang (2015) established a functional central limit theorem for a randomization of process $R_n$:
each indicator is multiplied independently by a
random variable taking values in $\pm 1$ with equal probabilities. The limiting Gaussian process is a sum of independent 
self-similar processes in this case.

Now we formulate the main result of the paper.

\begin{Theorem}

{\bf (i)} 
Let $\theta\in(0, 1)$ and $\nu\geq 1$ be an integer. 
Then process
$
\left(Y^{*}_{n,1}(t),\ldots,Y^{*}_{n,\nu}(t),
\ 0 \leq t \leq 1 \right)
$
converges weakly in the uniform metrics in  $D([0,1]^{\nu})$ 
to $\nu$-dimensional Gaussian process with zero expectation and covariance function 
$(c^{*}_{ij}(\tau,t))_{i,j=1}^{\nu}$: for $\tau\leq t$, $i,j \in \{ 1, \ldots, \nu\}$ (taking $0^0=1$)
\[
c^{*}_{ij}(\tau,t)=\left\{
\begin{array}{cr}
 \sum\limits_{s=0}^{i-1} 
\sum\limits_{m=0}^{j-s-1} \frac{\tau^s (t-\tau)^m K_{m+s,\theta}}{t^{m+s-\theta}s!m!}  
-
 \sum\limits_{s=0}^{i-1} 
\sum\limits_{m=0}^{j-1} \frac{\tau^s t^m K_{m+s,\theta}}{(t+\tau)^{m+s-\theta} s!m!}  , & i < j;\\
t^{\theta} \sum\limits_{m=0}^{j-1} \frac{K_{m,\theta}}{m!}  
-
 \sum\limits_{s=0}^{i-1} 
\sum\limits_{m=0}^{j-1} \frac{\tau^s t^m K_{m+s,\theta}}{(t+\tau)^{m+s-\theta} s!m!}, & i\ge j;\\
\end{array}
\right.
 \]
  $c^*_{ij}(\tau,t)=c^*_{ji}(t,\tau)$.

{\bf (ii)}
Let $\theta=1$. Then process 
$\left( \frac{R_{[nt]}-{\bf E} R_{[nt]}}{(nL^{*}(n))^{1/2}}, \ 0\le t \le 1
\right)$ 
converges weakly in the uniform metrics in $D(0,1)$ 
to a standard Wiener process.
\end{Theorem}

The limiting $\nu$-dimensional Gaussian process is self-similar with Hurst parameter $H=\theta/2<1/2$. 
Its first component coincides in distribution with the first component of the limiting process in Theorem~1 in Durieu and Wang (2015).  
The above Karlin's theorems  are partial 
 cases of Theorem 3 due to
$
c_{ij}=c^{*}_{ij}(1,1)-c^{*}_{i+1,j}(1,1)-c^{*}_{i,j+1}(1,1)+c^{*}_{i+1,j+1}(1,1).
$
In Section 2 we give a proof of Theorem~3. 

\section{Proof of Theorem~3}
\begin{Lemma} 
{\bf (i)} If $\theta \in (0,1)$ then there exist $n_0\ge 1$, $C(\theta)<\infty$ such that 
$
\frac{{\bf E} R_{\Pi(n\delta)}}{\alpha(n)} 
\le
C(\theta) \delta^{\theta/2}
$
for any $\delta\in[0,1], \ n\ge n_0$. If  $\theta=1$ then the same holds with  $n L^*(n)$ instead of $\alpha(n)$.
\\
{\bf (ii)}
Let $\tau\le t$, then ${\bf E} (R^*_{\Pi(t),k}-R^*_{\Pi(\tau),k})\le {\bf E} R_{\Pi(t-\tau)}$, $k\ge1$. 
\\
{\bf (iii)}
For any pair $\varepsilon$,
  $\delta\in(0,1)$ there exists $N=N(\varepsilon,\delta)$ such that for any $n\ge N$,
$
{\bf P}(\forall t\in[0,1] \  \ \exists \tau: |\tau-t|\le \delta,  
\ \Pi(n\tau)= [nt]) \stackrel{def}{=} {\bf P} (A(n))\ge1-\varepsilon.
$
\end{Lemma}

{\it Proof.} {\bf (i)} Let $\theta\in(0,1)$.
From Karamata representation  (Theorem 2.1, Appendix 6, inequality (A6.2.10) in Borovkov (2013)) exists  
$C_1(\theta)>0$ such that for all $x$ and $\delta\in(0,1]$ under condition $x\delta\ge C_1(\theta)$ 
there  is inequality $\frac{L(x\delta)}{L(x)}\le 2 \delta^{-\theta/2}$. 
 As
$\lim\limits_{x\to\infty}\frac{{\bf E} R_{\Pi(x)}}{\alpha(x)}=\Gamma(1-\theta)$ (Theorem 1 in Karlin (1967)), 
there exists $C_2(\theta)<\infty$ such that ${\bf E} R_{\Pi(x)}\le C_2(\theta) \alpha(x)$ as $x\ge x_0$ for some $x_0>1$. 

Let $n\delta>\max\{C_1(\theta),x_0\}$, then 
$
\frac{{\bf E} R_{\Pi(n\delta)}}{\alpha(n)} \le C_2(\theta)
\frac{(n\delta)^{\theta} L(n \delta)}{n^{\theta} L(n)}
\le
2 C_2(\theta) \delta^{\theta/2}.
$

If $n\delta \le \max\{C_1(\theta),x_0\}$ then
$
\frac{{\bf E} R_{\Pi(n\delta)}}{\alpha(n)} 
\leq
\frac{{\bf E} {\Pi(n\delta)}}{\alpha(n)} 
=\frac{n\delta}{n^{\theta}L(n)}.
$
Let $n_0$ such that for any $n\ge n_0$ we have $n^\theta L(n)\ge n^{\theta/2}$ then
\[
\frac{n\delta}{n^{\theta}L(n)}\le
\frac{n\delta}{n^{\theta/2}}
=(n\delta)^{1-\theta/2}\delta^{\theta/2}
\le (\max\{C_1(\theta),x_0\})^{1-\theta/2}\delta^{\theta/2}.
\]

If $\theta=1$ then we change $\alpha(n)$ to $n L^*(n)$, $L(n)$ to $L^*(n)$ and repeat the proof.
\\
{\bf (ii)}
$
{\bf E} (R^*_{\Pi(t),k}-R^*_{\Pi(\tau),k})
 = \sum\limits_{i=1}^{\infty} \sum\limits_{j=0}^{k-1} {\bf P}(\Pi_i(\tau)=j) {\bf P} (\Pi_i(t)-\Pi_i(\tau)\ge k-j) 
$
$
\le \sum\limits_{i=1}^{\infty}{\bf P} (\Pi_i(t-\tau)\ge1)={\bf E}R_{\Pi(t-\tau)}.
$
\\
{\bf (iii)} 
Let define $\Pi(x)=0$ for $x<0$. From monotonicity of Poisson process, it is enough to prove that  
for any pair $\varepsilon$,
  $\delta\in(0,1)$ there exists $N=N(\varepsilon,\delta)$ such that for any $n\ge N$,
$
{\bf P}(\forall t\in[0,1], \  \  \Pi(n(t-\delta))\le [nt]\le \Pi(n(t+\delta))) \ge 1-\varepsilon.
$

From Strong Law of Large Numbers (SLLN), for any  $\varepsilon, \ \delta \in(0,1)$ there exists 
$T_0=T_0(\varepsilon, \delta)$ such that
$
{\bf P}\left(\forall T\ge T_0, \ \frac{\Pi(T)}{T}\in\left[1-\frac{\delta}{4},1+\frac{\delta}{4}\right]\right)
\ge 1-\varepsilon.
$

Let $N= \frac{2}{\delta}\max(T_0,2)$, then $n(t+\delta)\ge n \delta\ge 2T_0$.
Then with probability not less then  $1-\varepsilon$ we have: for all $t \in [0,1]$
$$
\Pi(n(t+\delta))\ge n(t+\delta)\left(1-\frac\delta4\right)> n\left(t+\frac\delta2\right)>[nt].
$$

So, we need to prove only that
${\bf P}(\forall t\in[0,1], \  \ \Pi(n(t-\delta))\le [nt])\ge 1 - \varepsilon$.

If $t\in[0,\delta]$ then  $\Pi(n(t-\delta))=0\le [nt]$ a.s.

If $t\in[\delta,1]$ then $\Pi(n(t-\delta))\le \Pi(n(t-\frac\delta2))$ a.s., and
$n\left(t-\frac\delta2\right)\ge \frac{n\delta}2\ge T_0$, and
with probability not less then  $1-\varepsilon$ we have: for all $t \in [\delta,1]$
$$
\Pi(n(t-\delta))\le  n\left(t-\frac\delta2\right)\left(1+\frac\delta4\right)
\le nt-\frac{n\delta}4\le nt-1\le [nt].
$$

 Lemma 1 is proved.

{\bf Proof of Theorem 3}

{\bf Proof of (i). Step~1 (covariances)} 
Let $\tau\leq t$,
\begin{eqnarray*}
\widetilde{c}_{ij} (\tau,t)={\bf cov} (R^{*}_{\Pi(\tau),i}, R^{*}_{\Pi(t),j})
\\
=
\sum_{k=1}^{\infty} \bigg({\bf P} (\Pi_k(\tau)< i, \ \Pi_k(t)< j )-
{\bf P} (\Pi_k(\tau)< i) {\bf P} (\Pi_k(t)< j) \bigg).
\end{eqnarray*}

If $i< j$ then 
\[
\widetilde{c}_{ij} (\tau,t)=
\sum_{k=1}^{\infty} \sum_{s=0}^{i-1} \frac{(\tau p_k)^s}{s!}e^{-\tau p_k} 
\left(\sum_{m=0}^{j-s-1} \frac{((t-\tau)p_k)^{m}}{m!} e^{-(t-\tau)p_k} -
\sum_{m=0}^{j-1} \frac{(t p_k)^{m}}{m!} e^{-t p_k}\right) 
\]
\[
=
\int_0^{\infty} \sum_{s=0}^{i-1} \frac{\tau^s x^{-s}}{s!}e^{-\tau/x} 
\left(\sum_{m=0}^{j-s-1} \frac{(t-\tau)^m x^{-m}}{m!} e^{-(t-\tau)/x} -
\sum_{m=0}^{j-1} \frac{t^m x^{-m}}{m!} e^{-t/x}\right) d\alpha(x).
\]

We integrate by parts and divide into two integrals, then we make change $y=x/t$ in the first integral, 
$y=x/(t+\tau)$ in the second one:
\[
\widetilde{c}_{ij}(\tau,t) =
 \sum\limits_{s=0}^{i-1} 
\sum\limits_{m=0}^{j-s-1} \frac{\tau^s (t-\tau)^m t^{-m-s}}{s!m!}  
\int\limits_0^{\infty} ((m+s)y^{-m-s-1}-y^{-m-s-2})e^{-1/y} 
\alpha(ty) dy
\]
\[
-
 \sum\limits_{s=0}^{i-1} 
\sum\limits_{m=0}^{j-1} \frac{\tau^s t^m (t+\tau)^{-m-s}}{s!m!}  
\int\limits_0^{\infty} ((m+s)y^{-m-s-1}-y^{-m-s-2})e^{-1/y} 
\alpha((t+\tau)y) dy.
\]
Analogously for $i\geq j$,
\[
\widetilde{c}_{ij} (\tau,t)
=
\sum\limits_{m=0}^{j-1}
\bigg( \frac{1}{m!} \int\limits_0^{\infty}(m y^{-m-1}-y^{-m-2})e^{-1/y} \alpha(ty) dy 
\]
\[
-  \sum\limits_{s=0}^{i-1} 
\frac{t^m \tau^s (t+\tau)^{-m-s}}{s!m!}
\int\limits_0^{\infty} ((m+s)y^{-m-s-1}-y^{-m-s-2})e^{-1/y} 
\alpha((t+\tau)y) dy\bigg).
\]
For any integer $r\geq 0$,
\[
\int_0^{\infty} y^{-r-2}e^{-1/y} 
\alpha(ty) dy \sim 
\alpha(t) \int_0^{\infty} y^{\theta-r-2}e^{-1/y} 
dy =\alpha(t)\Gamma(r+1-\theta)
\]
as $t \to \infty$.
Note that $ \int_0^{\infty}(r  y^{\theta-r-1}- y^{\theta-r-2})e^{-1/y} dy=K_{r,\theta}$. 
So (because $\alpha(nt)/ \alpha(n) \to t^{\theta}$ as $n \to \infty$),
$
c^{*}_{ij}(\tau,t)=\lim_{n\to \infty} \widetilde{c}_{ij}(n\tau, nt)/ \alpha(n).
$

{\bf  Step~2 (convergence of finite-dimensional distributions)}
Analogously to proof of Theorem~1 in Dutko (1989) we have for any fixed
$m\geq 1$, $0<t_1<t_2<\ldots<t_m\leq 1$ triangle array of
$m\nu$-dimensional random vectors 
$
\left\{({\bf I}(\Pi_k(nt_j)\geq i)-{\bf P}(\Pi_k(nt_j)\geq i))\alpha^{-1/2}(n), \ i \leq \nu, \ j \leq m \right), 
k\leq n\}_{n\geq 1}
$
satisfies Lindeberg condition (see Borovkov (2009), Theorem 8.6.2, p.215).

{\bf Step~3 (relative compactness)}
Let for any $\tau_1\le \tau_2$,  
\[
R^{*}_{\Pi(\tau_2),k}-R^{*}_{\Pi(\tau_1),k}=\sum_{i=1}^{\infty} {\bf I}(\Pi_i(\tau_2)\ge k, \ \Pi_i(\tau_1)< k)
\stackrel{def}{=}\sum_{i=1}^{\infty} {\bf I}_{i}(\tau_1,\tau_2)=\sum_{i=1}^{\infty} {\bf I}_{i},
\] 
$ P_i=P_i(\tau_1,\tau_2)={\bf P}({\bf I}_{i})$. 
We will use designations ${\bf I}_i$ and corresponding $P_i$ for different values of $\tau_1<\tau_2$.

We need in a new process $Z_{n,k}^{**}(t)=\frac{R^{*}_{\Pi ([nt]),k}-{\bf E} R^{*}_{\Pi ([nt]),k}}{(\alpha(n))^{1/2}}$.

We (a) prove continuity of the limiting process; (b) prove that $Z_{n,k}^*$ and $Z_{n,k}^{**}$ are 'close';
(c) prove relative compactness of $Z_{n,k}^{**}$.

a) Let  $\tau_1=n t_1, \ \tau_2=n t_2$ for $t_1<t_2$, then
\[
{\bf E}(Z^*_{n,k}(t_2)-Z^*_{n,k}(t_1))^2= \sum\limits_{i=1}^{\infty} {\bf E} ({\bf I}_i-P_i)^2/\alpha(n)
\le \sum\limits_{i=1}^{\infty} P_i/\alpha(n)\le C(\theta) (t_2-t_1)^{\theta/2}.
\]
Above we used the fact that variance of an indicator is lesser than its expectation and Lemma 1(i,ii).
Using Step~1 and Theorem 1.4 in Adler (1990), we prove that the $k$-th component of the limiting Gaussian process is in $C(0,1)$ a.s.
So the limiting Gaussian process is in $C([0,1]^{\nu})$ a.s. weak convergence in Skorokhod topology implies the same in the uniform topology.

b)
As 
$R^*_{\Pi(nt),k}-R^*_{\Pi([nt]),k}\le \Pi([nt]+1)-\Pi([nt])$ a.s., and 
${\bf E}( \Pi([nt]+1)-\Pi([nt]))=1$
 we have  for any $\eta>0$
\[
{\bf P} (\sup\limits_{0\le t\le 1} |Z^*_{n,k}(t)-Z^{**}_{n,k}(t)|>\eta)\le
{\bf P} (\sup\limits_{0\le m \le n} (\Pi(m+1)-\Pi(m)+1)>\eta\sqrt{\alpha(n)})
\]
\[
\le \sum\limits_{m=0}^n {\bf E}e^{\Pi(m+1)-\Pi(m)+1}/ e^{\eta\sqrt{\alpha(n)}}=(n+1)e^{e-\eta\sqrt{\alpha(n)}}
\to 0
\] as $ n\to\infty$. So it is enough to prove relative compactness of $\{Z^{**}_{n,k}\}_{n\ge n_0}$.

c)
 Let $t_2-t_1\ge \frac{1}{2n}$, then 
$[nt_2]-[nt_1]\le n(t_2-t_1)+1 \le 3n(t_2-t_1)$.
Let $\gamma=[16/\theta]+1$, and $\tau_1=[n t_1], \ \tau_2=[n t_2]$. 
Using independence of terms and  Rosenthal inequality, we have for all $n\ge n_0$ (where $n_0$ is from Lemma 1 (i))
\[
{\bf E}|Z^{**}_{n,k}(t_2)-Z^{**}_{n,k}(t_1)|^{\gamma}\leq \frac{c(\gamma)}{(\alpha(n))^{\gamma/2}}
\left(
\sum\limits_{i=1}^{\infty} {\bf E} | {\bf I}_i-P_i|^{\gamma}+
\left( \sum\limits_{i=1}^{\infty} {\bf E} ( {\bf I}_i-P_i)^2\right)^{\gamma/2}\right)
\]
\[
\le
 \frac{c(\gamma)}{(\alpha(n))^{\gamma/2}}
\left(\sum\limits_{i=1}^{\infty} P_i + \left(
\sum\limits_{i=1}^{\infty} P_i \right)^{\gamma/2}\right)
\]
\[
\le
\frac{c(\gamma)}{(\alpha(n))^{\gamma/2}}\left(24n^4(t_2-t_1)^4+
({\bf E}R_{\Pi(2n(t_2-t_1))})^{\gamma/2}\right)
\le
\widetilde{C}(\theta)(t_2-t_1)^4.
\]

Here $c(\gamma)$ and $\widetilde{C}(\theta)$  depend on its argument only. 
Above we used the fact that variance of an indicator is lesser than its expectation,
inequality $\sum_{i} P_i \leq {\bf E} (\Pi([nt_2])-\Pi([nt_1]))\le 3n(t_2-t_1)\le 24n^4(t_2-t_1)^4$, and Lemma 1({i},{ii}).

Let $0\le t_2-t_1<1/n$,  then $[nt_1]=[nt]$ or $[nt_2]=[nt]$ for any $t\in[t_1, t_2]$.  
So 
$
Q\stackrel{def}{=} {\bf E} 
(|{Z}^{**}_{n,k}(t)-{Z}^{**}_{n,k}(t_1)|^{\gamma/2}
|{Z}^{**}_{n,k}(t_2)-{Z}^{**}_{n,k}(t)|^{\gamma/2})=0\le (t_2-t_1)^2.
$

Let $t_2-t_1\ge 1/n$,  then there are 3 possible cases:

1)  $t_2-t\ge \frac1{2n}$, $t-t_1\ge \frac1{2n}$, then 
from Cauchy-Bunyakovsky Inequality,  $Q\le \widetilde{C}(\theta) (t_2-t)^2(t-t_1)^2 \le \widetilde{C}(\theta)(t_2-t_1)^2$;

2) $t_2-t\ge \frac1{2n}$, $t-t_1< \frac1{2n}$, then 
from Cauchy-Bunyakovsky Inequality,  
\[
Q\le \sqrt{\widetilde{C}(\theta) (t_2-t)^4 {\bf E} \left(\frac{\Pi(1)+1}{\sqrt{\alpha(n)}}\right)^{\gamma}}\le \widehat{C}(\theta)(t_2-t_1)^2;
\]

3) $t_2-t< \frac1{2n}$, $t-t_1\ge \frac1{2n}$, symmetric to  case 2.

So we have (see Billingsley (1999), Theorem~13.5) density of $k$-th component and therefore density of all the vector.

{\bf Step~4 (approximation of the original process)}
From the relative compactness of distributions of processes $\{ Z^*_{n,k} \}_{n\ge n_0, k \ge 1}$ 
we get that for every pair $\varepsilon>0$,  $\eta>0$ there exist $\delta\in(0,1)$
and $N_1=N_1(\varepsilon,\eta)$ such that for all $n\ge N_1$
\[
 {\bf P}(\sup\limits_{|t-\tau| \leq \delta} \left|Z^*_{n,k}(\tau)-Z^*_{n,k}(t)\right|\ge \eta) \le \varepsilon.
\]
Then (as ${\bf P }(Y^*_{n,k}(t)=Z^*_{n,k}(\tau)|\Pi(n\tau)=[nt])=1$) we have for all $n\ge \max(N,N_1)$, where $N$ is from
Lemma 1 (iii),
\[
{\bf P}(\sup\limits_{0 \leq t \leq 1} \left|Y^*_{n,k}(t)-Z^*_{n,k}(t)\right|\ge \eta)\le 
{\bf P}(\sup\limits_{0 \leq t \leq 1} \left|Y^*_{n,k}(t)-Z^*_{n,k}(t)\right|\ge \eta, A(n))+\varepsilon
\]
\[
\le {\bf P}(\sup\limits_{|t-\tau| \leq \delta} \left|Z^*_{n,k}(\tau)-Z^*_{n,k}(t)\right|\ge \eta)+\varepsilon\le 2 \varepsilon.
\]
So (i) is proved.

{\bf Proof of (ii)}
Analogously to {\bf Step~1} for  $\tau<t$  
\[
\lim_{n\to \infty} \frac{{\bf cov} (R_{\Pi(n\tau)}, R_{\Pi(nt)} )}{n L^{*}(n)} = 
\lim_{n\to \infty} \frac{1}{n L^{*}(n)} \int_0^{\infty} e^{-nt/x}(1-e^{-n\tau/x})d \alpha (x) =\tau.
\]

Doing precisely the same as at {\bf Step~2--Step~4} (using slow variation of $L^{*}(x)$) we prove (ii).

{ \it The Theorem is proved.}

{\bf Acknowledgements}

The authors would like to thank Serguei Foss and an anonymous referee 
for their helpful and constructive comments and suggestions.

\bigskip

\footnotesize

{\sc Adler, R.J.}, 1990. An introduction to continuity, extrema, and related topics for general Gaussian processes,
Institute of Math. Stat., Hayward, California.


{\sc Barbour, A. D.}, 2009. Univariate approximations in the infinite occupancy
scheme. Alea 6, 415--433.

{\sc Barbour, A. D.,  Gnedin, A. V.}, 2009.
Small counts in the infinite occupancy scheme.
Electronic Journal of Probability, 
Vol. 14, Paper no. 13, 365--384.

{\sc Billingsley, P.}, 1999.
Convergence of Probability Measures, 2nd Edition, Wiley.

{\sc Borovkov, A. A.}, 2013.
Probability Theory,
Universitext.

{\sc Chebunin, M. G.}, 2014. Estimation of parameters of probabilistic models which is based on the number of different elements 
in a sample. Sib. Zh. Ind. Mat., 17:3, 135--147.


{\sc Durieu and Wang}, 2015. From infinite urn schemes to decompositions of self-similar Gaussian processes. Preprint.

{\sc Dutko, M.}, 1989.
Central limit theorems for infinite urn models, Ann. Probab. 17,
1255--1263.

{\sc Gnedin, A., Hansen, B., Pitman, J.}, 2007.
Notes on the occupancy problem with
infinitely many boxes: general
asymptotics and power laws.
Probability Surveys,
Vol. 4, 146--171.

{\sc Hwang, H.-K., Janson, S.}, 2008. Local Limit Theorems for Finite and Infinite 
Urn Models. The Annals of Probability, Vol. 36, No. 3,  992--1022.

{\sc Karlin, S.}, 1967. Central Limit Theorems for Certain Infinite Urn Schemes. 
Jounal of Mathematics and Mechanics, Vol. 17, No. 4,  373--401.

{\sc Key, E. S.}, 1992. Rare Numbers. Journal of Theoretical Probability,
Vol. 5, No. 2, 375--389.

{\sc Key, E. S.}, 1996.
Divergence rates for the number of rare numbers. Journal of Theoretical Probability, Volume
9, No. 2, 413--428.







{\sc Zakrevskaya, N. S.,  Kovalevskii, A. P.}, 2001. One-parameter probabilistic models of text statistics. 
Sib. Zh. Ind. Mat., 4:2, 142--153.

\end{document}